\documentstyle{amsppt}
\baselineskip18pt
\magnification=\magstep1
%\NoPageNumbers
%\NoRunningHeads
\pagewidth{30pc}
\pageheight{45pc}

\hyphenation{co-deter-min-ant co-deter-min-ants pa-ra-met-rised
pre-print pro-pa-gat-ing pro-pa-gate
fel-low-ship Cox-et-er dis-trib-ut-ive}
\def\leaderfill{\leaders\hbox to 1em{\hss.\hss}\hfill}
\def\A{{\Cal A}}

\def\H{{\Cal H}}
\def\L{{\Cal L}}

\def\e{{\varepsilon}}

\def\th{{\theta}}

\def\s{{\sigma}}

\def\P{{\widetilde P}}
\def\Q{{\widetilde Q}}

\def\b0{\text{\bf 0}}

\def\ra{{\ \longrightarrow \ }}

\def\lan{{\langle}}
\def\ran{{\rangle}}

\def\zed{{\Bbb Z}}

\def\boxit#1{\vbox{\hrule\hbox{\vrule \kern3pt
\vbox{\kern3pt\hbox{#1}\kern3pt}\kern3pt\vrule}\hrule}}
\def\rabbit{\vbox{\hbox{\kern0pt
\vbox{\kern0pt{\hbox{---}}\kern3.5pt}}}}

\topmatter

\title A projection property for Kazhdan--Lusztig bases
\endtitle

\author R.M. Green and J. Losonczy \endauthor
\vfill
\affil
Department of Mathematics and Statistics\\ Lancaster University\\
Lancaster LA1 4YF\\ England\\
{\it  E-mail:} r.m.green\@lancaster.ac.uk\\
\newline
Department of Mathematics\\
Long Island University\\
Brookville, NY  11548\\
USA\\
{\it  E-mail:} losonczy\@e-math.ams.org\\
\endaffil

\abstract
We compare the canonical basis for a generalized Temperley--Lieb
algebra of type $A$ or $B$ with the Kazhdan--Lusztig basis for 
the corresponding Hecke algebra.  
\endabstract

\thanks
The first author was supported in part by an award from the Nuffield
Foundation.
\endthanks

\endtopmatter

\centerline{\bf To appear in International Mathematics Research Notices}

\head Introduction \endhead

Generalized Temperley--Lieb algebras arise as certain quotients of
Hecke algebras associated to Coxeter systems, in the same way that the
ordinary Temperley--Lieb algebra can be realised as a quotient of
the Hecke algebra of type $A$ (see \cite{{\bf 9}}).  The finite 
dimensional generalized Temperley--Lieb algebras were classified 
by J. Graham \cite{{\bf 5}} into seven infinite families: 
types $A$, $B$, $D$, $E$, $F$, $H$ and $I$.

In \cite{{\bf 7}}, we showed that a generalized Temperley--Lieb algebra
arising from a Coxeter system of arbitrary type admits a canonical 
(more precisely, an IC) basis.
Such a basis is by definition unique; furthermore, it 
is analogous, in a manner which can be made precise, to G.\ Lusztig's 
canonical basis for the $\pm$-part of a quantized enveloping algebra 
and also to the Kazhdan--Lusztig basis for a Hecke algebra.  We 
determined the IC basis explicitly for the algebras of 
types $A$, $D$ and $E$; in each case the basis coincides with a 
previously familiar basis with convenient 
properties \cite{{\bf 7}, Theorem 3.6}.

In this paper, we examine the relationship between the IC
basis of a generalized Temperley--Lieb algebra and the corresponding
Kazhdan--Lusztig basis, showing in particular that in type $B$, the 
projection of a certain natural subset of the Kazhdan--Lusztig basis 
agrees with the IC basis.  Such a relationship is not obvious for a 
general Coxeter system. For example, the kernel of the relevant 
homomorphism need not be spanned by the Kazhdan--Lusztig basis elements 
which it contains.  We point out that some of the results in this 
paper lead to a new, concise and more general proof of the main result 
of \cite{{\bf 4}}.

\head 1. Preliminaries \endhead

\subhead 1.1 Generalized Temperley--Lieb algebras and 
IC bases \endsubhead

Let $X$ be a Coxeter graph, of arbitrary type,
and let $W(X)$ be the associated Coxeter group with distinguished
set of generating involutions $S(X)$.  Denote by $\H(X)$ the Hecke
algebra associated to $W(X)$.
(The reader is referred to \cite{{\bf 8}, \S7} for the basic theory of 
Hecke 
algebras arising from Coxeter systems.)  Let $\A$ be the ring of
Laurent polynomials, $\zed[v, v^{-1}]$.  The $\A$-algebra $\H(X)$ has 
a basis consisting of elements $T_w$, with $w$ ranging over $W(X)$, 
that satisfy 
$$T_s T_w = \cases
T_{sw} & \text{ if } \ell(sw) > \ell(w),\cr
q T_{sw} + (q-1) T_w & \text{ if } \ell(sw) < \ell(w),\cr
\endcases$$ where $\ell$ is the length function on the Coxeter group
$W(X)$, $w \in W(X)$, and $s \in S(X)$.
The parameter $q$ is equal to $v^2$.

Let $J(X)$ be the two-sided ideal of $\H(X)$ generated by the 
elements $$
\sum_{w \in \lan s, s' \ran} T_w,
$$ where $(s, s')$ runs over all pairs of elements of $S(X)$
that correspond to adjacent nodes in the Coxeter graph.  
(If the nodes corresponding to $(s, s')$ are connected by a
bond of infinite strength, then we omit the corresponding relation.)

\definition{Definition 1.1.1}
The generalized Temperley--Lieb algebra, $TL(X)$, is 
the quotient $\A$-algebra $\H(X)/J(X)$.  We denote the canonical 
epimorphism of algebras by $\th : \H(X) \ra TL(X)$.
\enddefinition

\definition{Definition 1.1.2}
A product $w_1w_2\cdots w_n$ of elements $w_i\in W(X)$ is called
{\it reduced} if $\ell(w_1w_2\cdots w_n)=\sum_i\ell(w_i)$.  We reserve
the terminology {\it reduced expression} for reduced products 
$w_1w_2\cdots w_n$ in which every $w_i \in S(X)$.

Call an element $w \in W(X)$ {\it complex} if it can be written 
as a reduced product $x_1 w_{ss'} x_2$, where $x_1, x_2 \in W(X)$ and
$w_{ss'}$ is the longest element of some rank 2 parabolic subgroup 
$\lan s, s'\ran$ such that $s$ and $s'$ do not commute.

Denote by $W_c(X)$ the set of all elements of $W(X)$
that are not complex.

Let $t_w$ denote the image of the basis element $T_w \in \H(X)$ in
the quotient $TL(X)$.
\enddefinition

\proclaim{Theorem 1.1.3 (Graham)}
The set $\{ t_w : w \in W_c \}$ 
is an $\A$-basis for the algebra $TL(X)$.  \endproclaim

\demo{Proof}
See \cite{{\bf 5}, Theorem 6.2}.
\qed\enddemo

The basis in Theorem 1.1.3 will be called the ``$t$-basis''.
The $t$-basis plays a r\^ole in the definition of other
important bases for the generalized Temperley--Lieb algebra.

\definition{Definition 1.1.4}
For each $s\in S(X)$, we define $b_s=v^{-1}t_s+v^{-1}t_e$.  If 
$w\in W_c(X)$ and $s_1s_2\cdots s_n$ is a reduced expression
for $w$, then we define $b_w = b_{s_1}b_{s_2}\cdots b_{s_n}$.

Note that $b_w$ does not depend on the choice of reduced
expression for $w$.

It is known that $\{b_w : w\in W_c\}$ is a basis for the $\A$-module
$TL(X)$.  We shall call it the {\it monomial basis}.
\enddefinition

We now recall a principal result of \cite{{\bf 7}}, which establishes
the canonical basis for $TL(X)$.  This basis is a direct analogue of 
the important Kazhdan--Lusztig basis of the Hecke algebra $\H(X)$.
  
Fix a Coxeter graph, $X$.
Let $(I,\leq)$ be the poset with $I=W_c(X)$ and with $\leq$ defined as
the restriction to $I$ of the Bruhat--Chevalley order on $W(X)$.

Let $\A^- = \zed[v^{-1}]$.  Let $\,\bar{\ }\,$ be the
involution on the ring $\A = {\Bbb Z}[v, v^{-1}]$ which satisfies 
$\bar{v} = v^{-1}$.

By \cite{{\bf 7}, Lemma 1.4}, the algebra $TL(X)$ has a $\zed$-linear 
automorphism of order $2$ which sends $v$ to $v^{-1}$ and $t_w$ to 
$t_{w^{-1}}^{-1}$.  We denote this map also by $\,\bar{\ }\,$.

Let $\L$ be the free $\A^-$-submodule of $TL(X)$ with basis 
$\{v^{-\ell(w)}t_w : w \in W_c\}$, and let $\pi : \L \ra \L/v^{-1}\L$ 
be the canonical projection.

\proclaim{Theorem 1.1.5}
There exists a unique basis $\{ c_w : w \in W_c\}$ for $\L$
such that $\overline{c_w} = c_w$ and $\pi(c_w) = \pi(v^{-\ell(w)}t_w)$
for all $w\in W_c$.
\endproclaim

\demo{Proof}
This is \cite{{\bf 7}, Theorem 2.3}.
\qed\enddemo

The basis $\{c_w : w \in W_c\}$ is called the {\it IC basis} (or the 
{\it canonical basis}) of $TL(X)$.  It depends on the 
$t$-basis, the involution $\,\bar{\ }\,$, and the lattice $\L$.  

\subhead 1.2 General questions \endsubhead

It was shown in \cite{{\bf 7}, Theorem 3.6} that the canonical basis of 
the generalized Temperley--Lieb algebra equals the monomial basis in 
types $A$, $D$ and $E$.  In this paper, we determine the
IC basis for type $B$.  Note that in any non-simply-laced 
case, the monomial basis does not equal
the canonical basis \cite{{\bf 7}, Remark 3.7 (1)}.
It would be interesting to have a description of the canonical basis 
in each of the finite dimensional types.

A natural related problem concerns the relationship
between the IC basis of a generalized Temperley--Lieb algebra
and the Kazhdan--Lusztig basis for the corresponding Hecke 
algebra (see \cite{{\bf 7}, Remark 2.4 (2)}).
Recall from \cite{{\bf 10}} that for each $w\in W$ there exists a
unique $C'_w \in \H$ such that $\overline{C'_w} = C'_w$, 
where $\,\bar{\ }\,$ is a certain $\zed$-linear automorphism of
$\H$, and  
$$
C'_w = \sum_{{x \in W} \atop {x \leq w}}\P_{x,w}(v^{-\ell(x)}T_x)
,$$ where $\P_{x,w}$ lies in $v^{-1} \A^-$ if $x < w$, and $\P_{w,w}
= 1$.

Let $\Cal C$ denote the set of $C'_w\in \H$ indexed by $w\in W_c$.
It is clear from \cite{{\bf 7}, Lemma 1.5} that the set $\th(\Cal C)$ 
is a basis for $TL(X)$.  

\definition{Definition 1.2.1}
We say that a Coxeter graph $X$ satisfies 
the {\it projection property} if $\th(\Cal C)$ equals the
canonical basis of $TL(X)$.
\enddefinition

We do not know of an example of a Coxeter graph which fails to
have the projection property.
The following two propositions provide useful sufficient conditions
for $\th(\Cal C)$ to equal the canonical basis.

\proclaim{Proposition 1.2.2}
The following are equivalent for a Coxeter graph $X$:

\item{\rm (i)}{ $\th(C'_w) \in \L$ for all $w \in W(X)$, where $\L$
is the lattice defined in \S1.1.}
\item{\rm (ii)}{ $\th(v^{-\ell(w)}T_w) \in \L$ for all $w \in W(X)$.}

If (i) or (ii) holds, then $\pi(\th(C'_w)) =
\pi(\th(v^{-\ell(w)}T_w))$ for all $w \in W$, so that $X$ satisfies 
the projection property.
\endproclaim

\demo{Proof}
For any $w \in W(X)$, we have $$
C'_w =\sum_{{x \in W(X)} \atop {x \leq w}}\P_{x, w}(v^{-\ell(x)} T_x)
,$$ where $\P_{x, w}$ lies in $v^{-1} \A^-$ if $x < w$, and $\P_{w, w}
= 1$.  If we assume (ii), we see that $\th(\P_{x, w} (v^{-\ell(x)}
T_x)) \in v^{-1}\L$ if $x < w$, and 
$\th(\P_{w, w} (v^{-\ell(w)} T_w)) = \th(v^{-\ell(w)} T_w) \in \L$.  
Thus $\th(C'_w) \in \L$ and $\pi(\th(C'_w)) =
\pi(\th(v^{-\ell(w)}T_w))$, which proves
that (ii) $\Rightarrow$ (i) and also that (ii) implies 
$\pi(\th(C'_w)) = \pi(\th(v^{-\ell(w)}T_w))$ for all $w \in W$.

Next, we show that (i) implies (ii).  
One sees from \cite{{\bf 10}} that for any $w \in W(X)$ we have $$
v^{-\ell(w)} T_w = \e_w \sum_{{x \in W} \atop {x \leq w}} \e_x
\Q_{x, w} C'_x,
$$ where $\e_x := (-1)^{\ell(x)}$, $\Q_{x, w}$ lies in $v^{-1} \A^-$
if $x < w$, and $\Q_{w, w} = 1$.  An argument similar to that of
the previous paragraph gives (i) $\Rightarrow$ (ii).

Finally, we verify that $X$ satisfies the projection property
if either (i) or (ii) holds.
It was shown in  \cite{{\bf 7}, Lemma 1.4} that the automorphism 
$\,\bar{\ }\,$ of $\H$ given in \cite{{\bf 10}} induces the 
automorphism $\,\bar{\ }\,$ of $TL(X)$ defined in \S1.1.  It
follows that $\overline{\th(C'_w)} = \th(C'_w)$.  
Since $\pi(\th(C'_w)) = \pi(\th(v^{-\ell(w)}T_w))$, Theorem 1.1.5 
shows that $\th(C'_w) = c_w$ if $w \in W_c$, which completes 
the proof.
\qed\enddemo

\proclaim{Proposition 1.2.3}
If $X$ is a Coxeter graph such that the kernel of the canonical map
$\th : \H(X) \ra TL(X)$ is spanned by the basis 
elements $C'_w$ which it contains, then $X$ has the projection 
property.
\endproclaim

\demo{Proof}
Suppose that the hypothesis is satisfied, and let $w \in W_c$.
Recall from the proof of Proposition 1.2.2 that $$
v^{-\ell(w)} T_w = \e_w \sum_{{x \in W} \atop {x \leq w}} \e_x
\Q_{x, w} C'_x
.$$  Applying $\th$ to both sides yields $$
\th(v^{-\ell(w)} T_w) = \e_w \sum_{{x \in W_c} \atop {x \leq w}} \e_x
\Q_{x, w} \th (C'_x)
.$$  This means that the change of basis matrix from
the basis $\{ \th(C'_w): w \in W_c\}$ to the basis 
$\{v^{-\ell(w)} t_w : w\in W_c\}$ (with respect to some total 
refinement of the partial order $\leq$) is upper triangular with ones 
on the diagonal, and all the entries above the diagonal lie in 
$v^{-1} \A^-$.  The inverse of this matrix has the same properties, 
meaning that $$
\th(C'_w)=\sum_{{x \in W_c} \atop {x \leq w}} \P_{x, w} 
(v^{-\ell(x)} t_x)
,$$ where $\P_{x, w} \in v^{-1} \A$ if $x < w$ and $\P_{w, w} = 1$.  
Thus $\pi(\th(C'_w)) = \pi(v^{-\ell(w)}t_w) = \pi(c_w)$, and  
since $\th(C'_w)$ is fixed by $\,\bar{\ }\,$, the proposition follows.
\qed\enddemo

We conclude this section with a conjecture 
(cf.\ \cite{{\bf 7}, Remark 2.4 (1)}).

\proclaim{Conjecture 1.2.4}
Let $X$ be an arbitrary Coxeter graph.  Then the structure constants 
of the canonical basis for $TL(X)$ lie in ${\Bbb N}[v, v^{-1}]$.
\endproclaim

The conjecture follows easily in types $A$, $D$ and $E$
from the results in \cite{{\bf 7}, \S3}.  

It is not difficult to see that when $X=I_2(m)$ (the dihedral case), 
the hypothesis of Proposition 1.2.3 is satisfied: the ideal 
$J(I_2(m))$ is spanned by the single Kazhdan--Lusztig basis element 
$C'_{w_0}$, where $w_0$ denotes the longest element of $W(I_2(m))$.
Moreover, it is known that the structure constants for the $C'$-basis 
of $\H(I_2(m))$ have positive coefficients.  Thus, one is able to 
deduce the positivity property for the IC basis of the generalized
Temperley--Lieb algebra $TL(I_2(m))$.

\head 2. Types $A$ and $B$ \endhead

Throughout this section, we assume that $X$ is a Coxeter graph 
of type $A$ or $B$. The algebra $TL(X)$ is then
generated by the monomial basis elements $b_s$, with $s$
ranging over $S(X)$, subject to the relations $b_s^2=q_cb_s$, 
$b_sb_{s'}=b_{s'}b_s$ if $ss'$ has order 2, 
$b_sb_{s'}b_s=b_s$ if $ss'$ has order 3, and 
$b_sb_{s'}b_sb_{s'}=2b_sb_{s'}$ if $ss'$ has order 4 
(see \cite{{\bf 6}, \S1}).  Here, $q_c$ is the Laurent polynomial 
$[2] = v + v^{-1}$.

Our goal is to prove that $X$ has the projection property.
This will be accomplished in Theorem 2.2.1.

\subhead 2.1 Reduced expressions \endsubhead 

Our first lemma describes a normal form reduced expression 
for elements in a Coxeter group of type $A$ or $B$.  

Let $\s_1,\s_2,\ldots ,\s_r$ denote the Coxeter generators of 
the Coxeter group $W(B_r)$, 
where $\s_1\s_2$ has order $4$ and $\s_i\s_{i+1}$ has order $3$ for
$i>1$. Define $W^{(r)}= 
\{w\in W(B_r) : 1\leq i<r\Rightarrow \ell(\s_iw)>\ell(w)\}.$
Then $W^{(r)}$ is a set of minimum length right
coset representatives for the subgroup $W(B_{r-1})$ of $W(B_r),$
and $\ell(xy)=\ell(x)+\ell(y)$ for all $x\in W(B_{r-1})$ 
and $y\in W^{(r)}$ (see \cite{{\bf 8}, \S5.12}).
Furthermore, it is known that the elements of $W^{(r)}$ are 
given by
$$\{
e,\, \s_r,\, \s_r\s_{r-1},\, \ldots ,\, \s_r\s_{r-1}\cdots 
\s_2\s_1,\, \s_r\s_{r-1}\cdots \s_2\s_1\s_2,\, \ldots, $$
$$\s_r\s_{r-1}\cdots \s_2\s_1\s_2\cdots \s_{r-1}\s_r 
\}.$$
(This can be established by working with
the signed permutation representation of $W(B_r)$.)
Observe that each element of $W^{(r)}$ has a unique 
reduced expression.  

\proclaim{Lemma 2.1.1} 
Let $w\in W(X)$.  There exists a reduced expression 
$s_1s_2\cdots s_n$ for $w$ such that for each $k$, either 
$s_k$ does not appear to the left of the $k$-th factor 
in $s_1s_2\cdots s_n$, or $s_k$ does not commute with $s_{k-1}$.
\endproclaim

\demo{Proof}
Every Coxeter system of type $B_r$ contains a Coxeter system of 
type $A_{r-1}$, the Coxeter generators of the latter being a subset 
of those of the former. Hence, it is enough to treat the case
$X=B_r$.

Considering the chain $\{e\}=W(B_0)\subseteq W(B_1)\subseteq \cdots
\subseteq W(B_r)=W$ of Coxeter groups, we see that 
any $w\in W$ has a unique reduced decomposition $w=w_1w_2\cdots w_r$, 
where each $w_i\in W^{(i)}$.  By substituting for each $w_i$ its 
unique reduced expression, we obtain a normal form reduced expression
for $w$ which satisfies the condition in the statement of the lemma. 
\qed \enddemo

Let $w\in W = W(X)$.  It is known~\cite{{\bf 1}, \S IV.1.5} that 
every reduced expression for $w$ can be transformed into any other
reduced expression for $w$ by performing a sequence of braid moves.  
Given this fact, one may characterize $W_c$ as the set of $w\in W$ 
such that every reduced expression for $w$ can be transformed into 
any other reduced expression for $w$ by a sequence of commutation 
moves (see~\cite{{\bf 11}, Proposition 1.1}).  

We define the {\it content} of $w\in W$ to be the set $c(w)$
of Coxeter generators $s\in S$ that appear in some 
(any) reduced expression for $w$.

The next result was stated without proof in \cite{{\bf 3}, \S 7.1}. 
The proof that we give is similar to that of \cite{{\bf 2}, Lemma 2}.

\proclaim{Lemma 2.1.2} 
Let $w\in W_c$ and $s\in S$ satisfy $ws\notin W_c$.  There
exists a unique $s'\in S$ such that
any reduced expression for $w$ can be parsed in one of the 
following two ways.

\item{\rm (i)}{$w=w_1sw_2s'w_3$, where $ss'$ has order 
$3$, and $s$ commutes with every member of 
$c(w_2)\cup c(w_3)$};
\item{\rm (ii)}{$w=w_1s'w_2sw_3s'w_4$, where 
$ss'$ has order $4$, $s$ commutes with every member of 
$c(w_3)\cup c(w_4)$, and $s'$ commutes with every member of 
$c(w_2)\cup c(w_3)$.}
\endproclaim

\demo{Proof}
A sequence of commutation moves applied to any reduced expression 
$x=s_1s_2\cdots s_n$ results in a reduced expression 
$x=s_{\tau(1)}s_{\tau(2)}\cdots s_{\tau(n)}$, 
where the $\tau(i)$-th factor of $s_1s_2\cdots s_n$ has been moved 
to the $i$-th position. 

We remark that if the generators $s_{\tau(i)}$ and 
$s_{\tau(j)}$ do not commute and $\tau(i)<\tau(j)$, then $i<j$.

Let $s_1s_2\cdots s_n$ be a reduced expression for $w$.  We shall
parse this reduced expression as in the statement of the lemma.

Define $s_{n+1}=s$.  Since $s_1s_2\cdots s_{n+1}\notin W_c$,
one can transform $s_1s_2\cdots s_{n+1}$ by a sequence 
of commutation moves into a reduced expression 
$s_{\tau(1)}s_{\tau(2)}\cdots s_{\tau(n+1)}$
which possesses a substring of type: (1) 
$s_{\tau(k)}s_{\tau(k+1)}s_{\tau(k+2)}$,
where $s_{\tau(k)}=s_{\tau(k+2)}$ and 
$(s_{\tau(k)}s_{\tau(k+1)})^3=1$; or type (2) $s_{\tau(k-1)}
s_{\tau(k)}s_{\tau(k+1)}s_{\tau(k+2)}$, where $s_{\tau(k-1)}=
s_{\tau(k+1)}$, $s_{\tau(k)}=s_{\tau(k+2)}$ and 
$(s_{\tau(k)}s_{\tau(k+1)})^4=1$.

We claim that in either case, $\tau(k+2)=n+1$.  By the
previous remark, we have $n+1\neq \tau(k), \tau(k+1)$, and also
$n+1\neq \tau(k-1)$ if we are considering a substring of type
(2). Hence, if $\tau(k+2)\neq n+1$, then we may commute $s_{n+1}$ 
to the end of $s_{\tau(1)}s_{\tau(2)}\cdots s_{\tau(n+1)}$ 
and thereby produce a reduced expression for $ws$ which contains a 
substring of type (1) or (2) among the first $n$ factors, 
contradicting $w\in W_c$.
Thus, we have $s=s_{\tau(k)}=s_{\tau(k+2)}$.

Again by the remark, we have $\tau(k)<\tau(k+1)<\tau(k+2)$,
and also $\tau(k-1)<\tau(k)$ if the substring is of type (2).
Suppose that the substring is of type (1).  Define $w_1 = 
s_1s_2\cdots s_{\tau(k)-1}$, $w_2 = s_{\tau(k)+1}s_{\tau(k)+2}
\cdots s_{\tau(k+1)-1}$, $s' = s_{\tau(k+1)}$ and 
$w_3 = s_{\tau(k+1)+1}s_{\tau(k+1)+2}\cdots s_n$.  Observe that in
order for the substring of type (1) to have occurred, it must 
be true that $s$ commutes with every member of $c(w_2)\cup c(w_3)$.
(To see why, consider the sequence of 
commutation moves required to transform the
subexpression $s w_2 s' w_3 s$ of $ws$ into an expression where the
occurrences of $s$, $s'$, $s$ are consecutive.)

If instead the substring is of type (2), then we define 
$w_1=s_1s_2\cdots s_{\tau(k-1)-1}$, $s'=s_{\tau(k-1)}$, 
$w_2 = s_{\tau(k-1)+1}s_{\tau(k-1)+2}\cdots s_{\tau(k)-1}$,
$w_3 = s_{\tau(k)+1}s_{\tau(k)+2}\cdots s_{\tau(k+1)-1}$ and 
$w_4 = s_{\tau(k+1)+1}s_{\tau(k+1)+2}\cdots s_n$.  As above,
in order for the substring of type (2) to have occurred, it 
must be the case that $s$ commutes with every member of 
$c(w_3)\cup c(w_4)$ 
and $s'$ commutes with every member of 
$c(w_2)\cup c(w_3)$. 

Regardless of the type of substring that arises, inspection of the
parsed form of the original reduced expression reveals that $s'$
is the rightmost factor which does not commute with $s$ in any 
reduced expression for $w$.  
Thus, $s'$ is uniquely determined by $w$ and $s$.
\qed \enddemo

Before presenting the next lemma, which furnishes some relevant
information concerning the structure constants of the monomial
basis, we remark that if $w=s_1s_2\cdots s_n$ is a reduced expression
and $1\leq i_1<i_2<\cdots <i_k\leq n$, then the product 
$b_{s_{i_1}}b_{s_{i_2}}\cdots b_{s_{i_k}}$ equals $aq_c^m b_x$, 
where $a$ and $m$ are nonnegative integers, and
$x\in W_c$ satisfies $\ell(x)\leq k$ and $x\leq w$.  This 
follows easily from the subexpression characterization of 
Bruhat--Chevalley order and the presentation of $TL(X)$ given at 
the beginning of this section.  We shall invoke this fact 
freely in the remainder of \S2.

We remind the reader that the elements $b_w\in TL(X)$ are always 
understood to be indexed by group elements $w\in W_c$ 
(see Definition 1.1.4). 

\proclaim{Lemma 2.1.3} 
Let $w\in W_c$ and $s\in S$.  Then $b_wb_s=aq_c^mb_{w'}$ for
some $w'\in W_c$ and some nonnegative integers $a$ and $m$, where
$m\leq 1$.  We have $\ell(w's)<\ell(w')$.  If 
$\ell(ws')<\ell(w)$ for some $s'\in S$ which does not commute
with $s$, then $m=0$.
\endproclaim

\demo{Proof}
We first address the last assertion in the statement.
Suppose that $w$ has a reduced expression ending
in a generator $s'$ which does not commute with $s$. Since
$w\in W_c$, this implies $\ell(ws)>\ell(w)$.  
If $ws\notin W_c$ then by Lemma 2.1.2, it is possible to write 
either $w=uss'$ (reduced), if $(ss')^3=1$, 
or $w=us'ss'$ (reduced), if $(ss')^4=1$.  In either 
case, $b_wb_s$ equals a positive integer multiple of some monomial 
basis element.

We now prove the rest of the lemma by induction on $\ell(w)$.  
The statement is true for $\ell(w)\leq 1$, so we may assume 
$\ell(w)>1$.

Suppose that $\ell(ws)<\ell(w)$.  We may write $w=us$ (reduced) 
for some $u\in W_c$.  We then have 
$b_wb_s=b_ub_sb_s=q_cb_ub_s=q_cb_w$, and the induction step follows.

Suppose instead that $\ell(ws)>\ell(w)$.  If $ws \in W_c$, the result
is obvious: set $m = 0$, $a = 1$ and $w' = ws$.  If $ws\notin W_c$
then we choose a reduced expression for $w$ as in the statement
of Lemma 2.1.1, and parse it according to Lemma 2.1.2. We treat 
only the case where the reduced expression may be parsed as
$w=w_1s'w_2sw_3s'w_4$, where $(ss')^4=1$, $s$ commutes with every 
member of $c(w_3)\cup c(w_4)$, and $s'$ commutes with every 
member of $c(w_2)\cup c(w_3)$. The reasoning
for the other case is similar. 

If $\ell(w_4)=0$, then $b_wb_s=
b_{w_1w_2w_3}b_{s'}b_sb_{s'}b_s=
2\,b_{w_1w_2w_3s's}$, and the induction step follows.

If $\ell(w_4)=1$, then we have $b_wb_s=
b_{w_1w_2w_3}b_{s'}b_sb_{s'}b_{w_4}b_{s}=
b_{w_1w_2w_3}b_{s'}b_sb_{s'}b_sb_{w_4}=2\,b_{w_1w_2w_3s's}b_{w_4}$.  
By induction, the last expression equals $aq_c^mb_{w'}$, 
where $m\leq 1$.  Also, since $b_wb_s=2b_{w_1w_2w_3s'}b_{w_4}b_s$,
induction gives $\ell(w's)<\ell(w')$.

If $\ell(w_4)>1$, then we parse $w_4=w_5s_1s_2$ ($s_1,s_2\in S$)
according to the chosen reduced expression for $w$, and then we 
apply the inductive hypothesis to 
$b_{ws_2}b_s$ (note that $\ell(ws_2)<\ell(w)$ and so $b_{w s_2}
b_{s_2} = b_w$), finding that it equals $aq_c^mb_{w'}$, where
$m\leq 1$ and $\ell(w's)<\ell(w')$.  Since our chosen reduced
expression for $w$ conforms to the condition in Lemma 2.1.1, either
$s_2\notin c(ws_2)$ or $s_2$ does not commute with $s_1$.  
In the first case, we have $$
b_wb_s = b_{ws_2}b_{s_2}b_s=b_{ws_2}b_sb_{s_2}=aq_c^mb_{w'}b_{s_2}=
aq_c^mb_{w's_2},
$$ where the last equality holds because $s_2 \not\in c(w')$.  
Since $s$ commutes with $s_2$, we have $\ell(w's_2s)<\ell(w's_2)$.

For the case where $s_1$ and $s_2$ do not commute, we
observe that $$
aq_c^m b_{w'} = b_{ws_2}b_s = 
b_{w_1w_2w_3}b_{s'}b_sb_{s'}b_{w_5}b_{s_1}b_s=
2\,b_{w_1w_2w_3s's}b_{w_5}b_{s_1}  
.$$  By induction,
$w'$ has a reduced expression ending in $s_1$.
Hence, by the argument given in the first paragraph,
if $w's_2\notin W_c$, then $b_{w'}b_{s_2}$ equals
$a'b_{w''}$ for some positive integer $a'$ and some $w''\in W_c$,
so that $b_wb_s =b_{ws_2}b_{s_2}b_s=b_{ws_2}b_sb_{s_2}=
aq_c^mb_{w'}b_{s_2}=aa'q_c^mb_{w''}$.  On the other hand, since 
$$b_wb_s = b_{w_1w_2w_3}b_{s'}b_sb_{s'}b_{w_4}b_s =
b_{w_1w_2w_3}b_{s'}b_sb_{s'}b_sb_{w_4}=2\,b_{w_1w_2w_3s'}b_{w_4}b_s 
,$$ induction gives
$\ell(w''s)<\ell(w'')$.

The induction step is complete.
\qed \enddemo

If a reduced expression satisfies the condition in the 
statement of the following lemma, then we say that it has 
the {\it deletion property}. 

\proclaim{Lemma 2.1.4} 
Let $w\in W$.  There exists a reduced expression 
$s_1s_2\cdots s_n$ for $w$ such that for any 
$1\leq i_1<i_2<\cdots <i_k\leq n$, we have 
$b_{s_{i_1}}b_{s_{i_2}}\cdots b_{s_{i_k}}=aq_c^mb_x$ for 
some $x\in W_c$ and some nonnegative 
integers $a$ and $m$, where $m\leq n-k$.
\endproclaim

\demo{Proof}
We shall prove by induction on $n=\ell(w)$ that if a reduced 
expression $s_1s_2\cdots s_n$ for $w$ satisfies the condition
in the statement of Lemma 2.1.1, then it has the deletion 
property.  When $n\leq 2$, this is evidently true, so we may assume 
$n>2$. Given $1\leq i_1<i_2<\cdots <i_k\leq n$, we consider 
two cases.

Case 1: $i_{k-1}<n-1$. We apply the inductive 
hypothesis to $s_1s_2\cdots s_{n-2}$ and the integers 
$1\leq i_1<i_2<\cdots <i_{k-1}\leq n-2$, thereby obtaining
$b_{s_{i_1}}b_{s_{i_2}}\cdots b_{s_{i_{k-1}}}=aq_c^mb_x$,
where $m\leq (n-2)-(k-1)$. By Lemma 2.1.3, 
we have $b_{s_{i_1}}b_{s_{i_2}}\cdots b_{s_{i_k}}=
aq_c^mb_xb_{s_{i_k}}=a'q_c^{m'}b_{x'}$, where 
$m'\leq m+1\leq n-k$.

Case 2: $i_{k-1}=n-1$. Here, we apply induction to
$s_1s_2\cdots s_{n-1}$ and the integers 
$1\leq i_1<i_2<\cdots <i_{k-1}=n-1$.  We obtain
$b_{s_{i_1}}b_{s_{i_2}}\cdots b_{s_{i_{k-1}}}=aq_c^mb_x$,
where $m\leq (n-1)-(k-1)$.  If $s_n\notin c(x)$, then 
the conclusion follows after multiplying both sides of the last 
equation on the right by $b_{s_n}$. If instead 
$s_n\in c(x)$, then because our reduced expression for $w$ satisfies 
the condition in Lemma 2.1.1,
$s_{n-1}$ and $s_n$ do not commute.  Furthermore, by Lemma 2.1.3, 
the element $x\in W_c$ has a reduced expression ending in 
$s_{i_{k-1}}=s_{n-1}$, hence by the last assertion of the same lemma, 
$b_xb_{s_n}$ equals a nonnegative integer multiple of some monomial 
basis element.
\qed \enddemo

We note that Lemma 2.1.4 may also be proved using the diagram calculus
for $TL(B_n)$ described in \cite{{\bf 6}, Theorem 4.1} in conjunction 
with Lemma 2.1.1.

\subhead 2.2 Canonical basis in types $A$ and $B$ \endsubhead

We are now prepared to state and prove the main result of this
paper, which can be viewed as a generalization of 
\cite{{\bf 4}, Theorem 3.8.2}.

\proclaim{Theorem 2.2.1} 
Let $X$ be a Coxeter graph of type $A$ or $B$.  Then the canonical 
basis of $TL(X)$ equals the image of the set of 
Kazhdan--Lusztig basis elements $C'_w\in \H(X)$ indexed by 
$w\in W_c(X)$. 
\endproclaim

\demo{Proof}
By Proposition 1.2.2 (ii), it suffices to show that for all 
$w\in W$, the element $v^{-\ell(w)}t_w$ lies in the lattice $\L$.

Every $w\in W$ has a reduced expression $s_1s_2\cdots s_n$ 
as in the statement of Lemma 2.1.4.  We have
$$\eqalign{
v^{-\ell(w)}t_w &= 
(v^{-1}t_{s_1})(v^{-1}t_{s_2})\cdots (v^{-1}t_{s_n}) \cr
 &= 
(b_{s_1}-v^{-1})(b_{s_2}-v^{-1})\cdots (b_{s_n}-v^{-1}). \cr
}$$
The last product expands to a sum of terms
$$
\pm v^{k-\ell(w)}b_{s_{i_1}}b_{s_{i_2}}\cdots b_{s_{i_k}},
$$
and each such term equals 
$$
a v^{k-\ell(w)}q_c^mb_x,
$$
where $a\in \zed$, $x\in W_c$ satisfies $x\leq w$, and 
$m\leq \ell(w)-k$ by Lemma 2.1.4.  The coefficient of $b_x$ in 
the expansion therefore lies in $\A^-$.  Also, one
sees that if $w\in W_c$, then the coefficient of $b_w$ equals 1.

We have shown that for any $w\in W$, the element $v^{-\ell(w)}t_w$
equals a linear combination of $b_x$ ($x\in W_c$ and $x\leq w$) 
with coefficients in $\A^-$; and if $w\in W_c$, then the
coefficient of $b_w$ equals 1. In particular, 
the change of basis matrix from the monomial basis 
$\{b_x : x\in W_c\}$ to the basis 
$\{v^{-\ell(x)}t_x : x\in W_c\}$ (with respect to some total 
refinement of $\leq$) is upper triangular with entries in $\A^-$
and with ones on the diagonal.  The inverse of this matrix must also  
have this property.  That is, for any $x\in W_c$, the element $b_x$ 
equals a linear combination of $v^{-\ell(y)}t_y$ 
($y\in W_c$ and $y\leq x$) with coefficients in 
$\A^-$.  Combining this with the first statement of the present 
paragraph, we conclude that every $v^{-\ell(w)}t_w$ lies in $\L$.
\qed \enddemo

\leftheadtext{}
\rightheadtext{}
%\vfill\eject

\end